\documentclass
[12pt,a4paper]{article}
\usepackage[cp1251]{inputenc}
\usepackage{amssymb, amsfonts, amsmath, latexsym}
\usepackage[russian]{babel}
\begin{document}

\Large 
\centerline{\bf  Ideals in $\mathcal{P} _G$  and $\beta G$ }\vspace{6 mm}

\normalsize\centerline{\bf  Igor Protasov and Ksenia Protasova}\vspace{6 mm}

{\bf Abstract.} For a discrete group $G$, we use the natural correspondence between ideals in the Boolean algebra $ \mathcal{P}_G$ of subsets of $G$ and closed subsets in the Stone-$\check{C}$ech compactifi-cation $\beta G$ as a right topological semigroup to introduce and  characterize some new ideals in $\beta G$. We show that if a group $G$ is either countable or Abelian then there are no closed ideals in $\beta G$ maximal in $G^*$, $G^* = \beta G \setminus G$, but this statement does not hold for the group $S_\kappa$ of all permutations of an infinite cardinal $\kappa$. We characterize the minimal closed ideal in $\beta G$ containing all idempotents of $G^*$. 

\vspace{6 mm}

2010 MSC:  22A15, 03E05

\vspace{3 mm}

Keywords: Stone-$\check{C}$ech compactification, Boolean algebra, ideal, filter, ultrafilter.
\vspace{6 mm}


\newcommand{\IN}{\mathbb N}

\section{Introduction}

We recall that a family $ \mathcal{I}$ of subsets of a set $X$ is an {\it ideal}  in the Boolean algebra $\mathcal{P} _G$ of all subsets of $G$ if $\emptyset \notin \mathcal{I}$ 
and $A\in \mathcal{I}$, $B\in \mathcal{I}$, $C\subseteq A$  imply $A\cup B\in \mathcal{I}$, $C\in \mathcal{I}$.  A family  $\varphi$ of subsets of $G$ is a {\it filter} if the family $\{ X\setminus  A: A \in \varphi \}$ is an ideal. A filter maximal by the inclusion is called  an {\it ultrafilter}.

For an infinite group $G$, an ideal $\mathcal{I}$ in $ \mathcal{P}_G$ is called {\it left (right) translation invariant} if $gA\in \mathcal{I}$ ($Ag\in \mathcal{I}$) for all $g\in G$, $A\in \mathcal{I}$. If $\mathcal{I}$ is left and right translation invariant then $\mathcal{I}$ is called {\it translation invariant}. Clearly, each left (right) translation invariant ideal of $G$ contains the ideal  $\mathcal{F} _G$ of all finite subsets of $G$. An ideal $\mathcal{I}$ in $\mathcal{P} _G$ is called a {\it group ideal} if  $\mathcal{F} _G\subseteq \mathcal{I}$ 
and if $A\in \mathcal{I}$, $B\in \mathcal{I}$  then $AB ^{-1}\in \mathcal{I}$. 

Now we endow $G$ with the discrete topology and identify the Stone-$\check{C}$ech compactifi-cation of $G$ with the set of all ultrafilters on $G$ and denote $G^{*}=\beta G\setminus G$, so $G^{*}$  is the set of all free ultrafilters on $G$. Then the family $ \{ \overline{A} : A \subseteq G \}$, where $\overline{A} = \{ p\in \beta G : A\in p \}$ forms the base for the topology of $\beta G$. Given a filter $\varphi$ on $G$, we denote $\overline{\varphi} = \cap \{ \overline{A} : A\in \varphi \}$, so $\varphi$  defines the closed subset   $\overline{\varphi}$ of   $\beta G$, and each non-empty closed subset $K$ of  $\beta G$  can be defined in this way: $K=\overline{\varphi} $  where $ \overline{\varphi} = \{ A\subseteq G:  K \subseteq   \overline{A} \}$.

We use the standard extension [4, Section 4.1] of the  multiplication on $G$ to the semigroup multiplication on  $\beta G$  such  that, for each $p\in \beta G$, the mapping $x\longmapsto  xp$, $x\in \beta G$  is continuous, and for each $g\in G$, the mapping, $x\longmapsto gx$, $x\in \beta G$ is continuous. Given  two ultrafilters $p, q\in \beta G$, we choose  $P\in p$ and, for each $x\in P$, pick $Q_x\in q$. Then $\bigcup  _{x\in P} \ xQ_x \in pq$  and the family of all these subsets forms the base of the product $pq$.
 \vskip 5pt

It follows directly from the definition of the multiplication in $\beta G$ that   $G^{*}$,   $ \overline{ G^{*} G^{*}}$ are ideals  in $\beta G$, and  $G^{*}$  is the unique  maximal closed ideal in  $\beta G$. By Theorem 4.44 from [4],  the closure  $ \overline{ K(\beta G)}$  of the minimal ideal    $K(G)$  of   $\beta G$  is an ideal, so $ \overline{ K(\beta G)}$ is the smallest closed ideal in $\beta G$.  For the structure of    $ \overline{ K(\beta G)}$  and some other ideals in  $\beta G$  see [4, Sections 4,6].

For an ideal $\mathcal{I}$ in $\mathcal{P}_G$, we put  
$$\mathcal{I}^{\wedge} = \{p\in  \beta G:  p\in G\setminus A  \text{  for each  }  A\in \mathcal{I} \}, $$
 and use the following observations: \vskip 5pt

\begin{itemize}
\item{} $\mathcal{I}$ is left translation invariant if and only if $\mathcal{I}^{\wedge}$    is a left ideal of the semigroup $\beta G$ ; \vskip 5pt

\item{}  $\mathcal{I}$ is right  translation invariant if and only if 
$(\mathcal{I}^{\wedge})G\subseteq  \mathcal{I}^{\wedge}$.

\end{itemize}
\vskip 5pt

We use also the inverse to $^\wedge$  mapping  $^\vee$.  For a closed  subset $K$ of $\beta G$, we take a filter $\varphi$  on $G$  such that $K=\overline{\varphi}$  and put 
$$K^{\vee} = \{G\setminus A : A\in    \varphi\}.$$

In   section 2, we use a classification of subsets of a group by their size to define some special ideals in
 $\mathcal{P}_G$.
 In section 3, we study ideals of  $\beta G$  between     $ \overline{ G^{*} G^{*}}$ 
and $G^{*}$. In section 4, we study ideals between $ \overline{ K(\beta G)}$ 
 and $ \overline{ G^{*} G^{*}}$  and characterize the minimal closed ideal in $\beta G$
  containing all idempotents  of  $G^{*}$.

\section{Diversity of subsets of a group}

In what follows, all group are supposed to be infinite. Let $G$  be a group with the identity $e$. We say that a subset $A$ of $G$ is \vskip 5pt

 \begin{itemize}
\item{} {\it large} if  $G=FA$  for some $F\in \mathcal{F}_G$;\vskip 5pt

\item{} {\it small} if  $L\setminus A$ is  large for every large subset $L$;\vskip 5pt

\item{} {\it thin} if  $gA\cap A$ is  finite for  each $g\in G\setminus \{e\}$;\vskip 5pt

\item{} {\it $n$-thin}, $n\in {\mathbb N}$
 if, for every   distinct elements   $g_0 , \dots , g_n \in  G$,
the set  $g_0  A \cap \dots \cap g_n  A$    is  finite; \vskip 5pt

\item{} {\it sparse} if,  for every   infinite subset   $X$ of $G$, there exists a  finite subset $F\subset X$  such that $\bigcap _{g\in F}  gA$ is  finite. \vskip 5pt
\end{itemize}

All above definitions can be unified with usage the following notion [16]. Given a subset $A$ of a group $G$ and an ultrafilter $p\in G^*$, we define a $p$-{\it companion} of $A$  by
$$ \Delta _p (A)= A^* \cap Gp=\{gp: g\in G,  \   A\in gp  \}.  $$
\vskip 5pt

Then the following statement hold [16]:
\vskip 5pt

 \begin{itemize}
\item{} $A$ is large if  and only if  $\Delta _p (A) \neq  \emptyset  $ for each    $p\in G^*$;\vskip 5pt

\item{} $A$ is  small if   and only if,  
for every  $p\in   G^*$  and every    $F\in  \mathcal{F}_G$,  we have $\Delta _p (F A) \neq  Gp  $;\vskip 5pt

\item{} $A$ is thin if and only if, $\Delta _p (A) \leq 1 $  for every $p\in G^*$;\vskip 5pt

\item{} $A$ is $n$-thin if and only if,  $\Delta _p (A) \leq n $   for every $p\in G^*$; \vskip 5pt

\item{} $A$ is sparse if and only if, $\Delta _p (A) $ is finite for each $p\in G^*$.\vskip 5pt

\end{itemize}

 Following [1], we say that a subset $A$  of $G$ is {\it scattered} if, for every infinite subset $X$ of $A$, there is $p\in X^* $  such that $\Delta _p (X)$ is finite. Equivalently [1, Theorem 1],  $A$ is scattered if each subset  $\Delta _p (A)$ is discrete in $G ^* $. 

We denote by $Sm_G$, $Sc_G$, $Sp_G$ the families of all small, scattered and sparse subsets of a group $G$. These families are translation invariant ideals in $\mathcal{P}_G$
 (see [16, Proposition 1 ]),  and for every group $G$, the following inclusions are strict [16,  Proposition 12]   
$$ Sp_G \subset  Sc_G  \subset  Sm_G . $$

We say that a subset $A$  of $G$ is {\it finitely thin} if $A$  is $n$-thin for some $n\in \IN$.
The family $ FT_G$ of  all  finitely  thin subsets of $G$  is a translation invariant ideal in   
$\mathcal{P}_G$   which
 contains the ideal $<T_G>$ generated by the family of all thin subsets of $G$.
 By  [6, Theorem 1.2] and [14, Theorem 3], if $G$ is either countable or Abelian and $|G| < \aleph _\omega$ then 
$FT_G = <T_G>$. By [14, Example 3], there exists a group $G$ of cardinality  $\aleph _\omega$
such that $<T_G>\subset FT_G $. 

Clearly, $FT_G \subseteq Sp_G $.  In the next section, we show that
$FT_G \subset Sp_G $  for every group  $G$
\vskip 5pt

{\bf Theorem 2.1.}
{\it For every group $G$, we have $Sm_G ^\wedge
= \overline{ K(\beta G)}$. }\vskip 5pt

This is Theorem 4.40 from [4]  in the form given in [10, Theorem 12.5].

\vskip 5pt
{\bf Theorem 2.2.}
{\it  For every group $G$, the following statements hold: 
\vskip 5pt

$(i)$    $Sp_G ^\wedge= \overline{ G ^*G ^*}$;
\vskip 5pt

$(ii)$ for a subset $A$ of $G$, $\overline{ G ^*G ^*} \subset \overline{ A}$   if and only if, for any infinite subsets $X, Y$ of $G$, there exist $x\in X$, $y\in Y$ such that  $xy\in A$,  $yx\in A$. }
\vskip 5pt

The statement $(i)$ is Theorem 10  from [2],  $(ii)$  is a recent result [11]. 

For more delicate classifications  of subsets of groups  and $G$-spaces see [5], [9], [15].


\section{Between $\overline{ G ^*G ^*}$  and  $G ^*$}

{\bf Theorem 3.1.}
{\it  For every group $G$, the following statements hold: 
\vskip 5pt

$(i)$   if  $\mathcal{I}$ is a  left  translation invariant ideal in $\mathcal{P}_G$  and $\mathcal{I} \neq \mathcal{F}_G$ then there exists a left translation invariant ideal  $\mathcal{J}$ in  $\mathcal{P}_G$ such  that  $\mathcal{F}_G \subset \mathcal{J} \subset  \mathcal{I} $ and   $\mathcal{J} \subset  Sp_G $;
\vskip 5pt

$(ii)$  if  $\mathcal{I}$ is a  right  translation invariant ideal in $\mathcal{P}_G$  and $\mathcal{I} \neq \mathcal{F}_G$ then there exists a right  translation invariant  $\mathcal{J}$ in  $\mathcal{P}_G$ such  that  $\mathcal{F}_G \subset \mathcal{J} \subset  \mathcal{I} $;
\vskip 5pt

$(iii)$  if $G$  is  either countable or Abelian and 
$\mathcal{I} $ is a  translation invariant ideal in $\mathcal{P}_G$ such  that
 $\mathcal{I} \neq \mathcal{F}_G$ then there exists a   translation invariant   ideal  $\mathcal{J}$ in  $\mathcal{P}_G$ such  that  $\mathcal{F}_G \subset \mathcal{J} \subset  \mathcal{I} $ and   $\mathcal{J} \subset  Sp_G $;
\vskip 5pt

\vskip 5pt

Proof.} We use the following auxiliary statement [8. Example 3]:
\vskip 5pt

$(*)$  if a countable group $\Gamma $ acts on a set $X$  then, for every infinite subset $A$  of $X$, there exists a countable subset $T\subset A$  such that the set $$\{ x\in T: gx\neq x, \  gx\in T ,  \    g\in G   \}$$ is finite.
\vskip 5pt

$(i)$   We suppose that  $G$ is countable, put  $\Gamma =G $,     $X=G$
 and consider the action of  $G$ on $G$ by the left shifts. We take an infinite  subset  $A\in  \mathcal{I}$
  and apply $(*) $  to choose a countable thin subset $T\subseteq A$. 
We partition $T$ into two infinite subsets $T=B\cup C$  and denote 
$$\mathcal{J}= \{ Z\subseteq G : \  Z\neq \emptyset,  \  Z\subset FB   \  \text{  for some }
\   F\in  \mathcal{F}_G \}.$$

Clearly,  $\mathcal{J}$ is a left translation invariant ideal and  $\mathcal{J}  \neq   \mathcal{F}_G$. 
Since  $gT\cap T$  is finite for every $g\in G\setminus \{e\}$,  we have $C\notin \mathcal{I} $.

Hence, $\mathcal{F}_G  \subset  \mathcal{J}  \subset  \mathcal{I} $. 
By the choice of $T$,  each subset $Y\in \mathcal{J} $
 is a finite union of thin subsets, so   $\mathcal{J}\subset  Sp_G $.
 
If $G$  is an arbitrary infinite group then we take a countable subset    $A\in \mathcal{I} $,
 consider the subgroup $H$ of $G$ generated by $A$  and denote by $\mathcal{I} _H$
 the restriction of  $\mathcal{I} $  to $H$,    $\mathcal{I} _H =\{ Y\cap H: Y\cap H \in  \mathcal{I}  \}  $.
 By above paragraph, there exists a left invariant ideal  
$\mathcal{J} ^\prime $ in $\mathcal{P} _H  $ such that 
$\mathcal{F} _H\subset  \mathcal{J} ^\prime \subset  \mathcal{I}_H$, 
$\mathcal{J} ^\prime \subset   Sp_H$.  
Then we put  $\mathcal{J}= \{Y \subseteq G : Y\neq \emptyset ,  \  Y\subseteq FZ,  \  F\in  \mathcal{F} _G , \  Z\in   \mathcal{J} ^\prime \}$.
  
\vskip 10pt

$(ii)$ We repeat the proof of $(i)$  with the action of $G$ on $G$  by the right shifts. 
\vskip 10pt

$(iii)$  If $G$  is countable then we put $\Gamma = G\times  G$, consider the action of $\Gamma$  on $G$ defined by  $(g,h)x= g^{-1} xh $
 and repeat the proof of $(i)$  in the countable case. If $G$  is Abelian then we apply $(i)$ directly.   $ \ \  \ \ \Box$

\vskip 10pt
{\bf Theorem 3.2.}
{\it  For every group $G$, the following statements hold: 
\vskip 5pt

$(i)$   if  $L$ is a  closed left ideal in   $\beta G$ such that $L\subset G^*$ then there exists a closed left ideal 
$L ^\prime$  of    $\beta G$   such that $L\subset L^\prime  \subset G^* $,   $\overline{ G ^*G ^*}\subset  L ^\prime$;

\vskip 5pt

$(ii)$   if  $R$ is a  closed subset of  $G^*$ such that $R\neq G^*$ and  $RG\subseteq  R$
 then there exists a closed subset
$R ^\prime$  of    $G^*$   such that
 $R\subset R^\prime  \subset G^* $,   $ R ^\prime G\subseteq  R $;

\vskip 5pt

$(iii)$   if  $G$ is  either countable or Abelian and  $I$  is 
a  closed 
ideal in   $\beta G$ such that 
$I\subset G^*$ then there exists a closed ideal 
$I ^\prime$  in    $\beta G$   such that $I\subset I^\prime  \subset G^* $,   
$\overline{ G ^*G ^*}\subset  I $.
\vskip 15pt

Proof.}  $(i)$ We put $\mathcal{I}= L^{\vee}$, apply Theorem 3.2 (i)  and set  $L^{\prime}=\mathcal{J}^{\vee}$. Then $L^{\prime}$ is a left ideal in $\beta G$  and 
$L \subset L^\prime  \subset G^* $. 
Since $\mathcal{J} \subset Sp_G$,  by Theorem 2.2, we have  $\overline{ G ^*G ^*}\subset   L^\prime $.

\vskip 10pt

$(ii)$ We put $\mathcal{I}= R^{\vee}$, and note that $\mathcal{I}$ is right translation   invariant. We apply Theorem 3.2(ii) and set $ R^{\prime}  = \mathcal{J}^{\wedge}$. 
\vskip 10pt 

$(iii)$ We put $\mathcal{I}= I^{\vee}$,
apply Theorem 3.2 (iii)  and set  
$I^{\prime}=\mathcal{J}^{\wedge}$. Then $I^{\prime}$ is a left ideal in $\beta G$  and 
$ I^\prime G \subseteq  I^\prime  $. 
Since $\mathcal{J} \subset Sp_G$, 
 we have  $I^\prime  G    \subseteq   I^\prime $ so  $I^\prime $  is a right ideal.

\vskip 15pt  

{\bf Remark 3.1.}  If  $\mathcal{I}$ is a group ideal in $\mathcal{P}_G$ then, by [13],  
 $\mathcal{I}^{\wedge}$  
is an ideal in $\beta G$. By [12, Theorem 4], if $G$ is either contable or Abelian and    $\mathcal{I}$  is a group ideal such that  $\mathcal{I}\neq \mathcal{F}_G$ then there exists a group ideal $\mathcal{J}$ in  $\mathcal{P}_G$  such that  $\mathcal{F}_G\subset  \mathcal{J} \subset  \mathcal{I}$. If $A$ is an infinite subset of $G$  then the subset 
$AA$ is not sparse  (put $X=A^{-1}$ in corresponding definition).  It follows that if $\mathcal{I}$ is a group ideal and   $\mathcal{I} \subseteq Sp_G$  then   $\mathcal{I}= \mathcal{F}_G$ .  

For a cardinal $\kappa$, $S_\kappa$ denotes the group of all permutations  of $\kappa$. 
\vskip 15pt

{\bf Theorem 3.3.}
{\it  For every infinite  cardinal $\kappa$,  there exists a closed ideal $I$ in $\beta S_\kappa$ such that 
\vskip 10pt 

$(i)$ $S_\kappa ^*  S_{\kappa ^*}\subset I  $; 
\vskip 10pt 

$(ii)$ if $M$ is a closed ideal in $\beta S_\kappa $ and $I\subseteq M\subseteq  G^*$  then either $M=I$ or $M= S_\kappa ^*$
\vskip 10pt

Proof.} 
We take an arbitrary closed subset $X= \{x_i : i<\omega \}$ of $\kappa$ and define a permutation $f_i$ of $\kappa$ by $f_i (x_{2i} ) = x_{2i+1} $,  $ \  f_i (x_{2i+1} ) = x_{2i} $ and  $f_i (x) = x$ for all $x\in \kappa \setminus \{x_{2i} ,   x_{2i+1} \}$. We put $T= \{f_i : i<\omega \}$ and denote by  $\mathcal{I}$  the smallest translation invariant ideal in  $\mathcal{P}_{S_\kappa}$
containing $T$. 

We note that $|gT \cap T|\leq 1$  for every  $g\in G\setminus \{e \}$.  Hence, $T$ is thin and $\mathcal{I}\subseteq Sp_{S_\kappa}$. 
To see that $\mathcal{I} \subset Sp_{S_\kappa}$, 
we observe that each element of   $\mathcal{I}$   is a countable subset of $S_\kappa$, but there are uncountable thin subsets of $S_\kappa$. 

We assume  that there is a translation invariant ideal $\mathcal{J}$ in $\mathcal{P}_{S_\kappa}$  such that 
$\mathcal{F}_G \subset \mathcal{J}  \subset \mathcal{I} $. 
Then there exists a countable subset $T_1$ of $T$ such that  $T_1\in \mathcal{J}$,  $ T\setminus  T_1 $   is infinite and  $ T\setminus  T_1  \notin \mathcal{I} $ .  We denote  $T_2= T\setminus  T_1$ and take a partition $\omega = W_1 \cup W_2$  such that  $T_1  = \{f_i : i\in W_1  \}$,    $T_2  = \{f_i : i\in W_2  \}$. We fix an arbitrary bijection $\varphi : W_1  \rightarrow  W_2  $ and define a permutation $h$  of $\kappa$ by the following rule.

 If $x\in \kappa\setminus X$  then $f(x)=x$.

If $x\in X$ then we take $i<\omega$   such that  $x\in \{ x_{2i} , x_{2i+1} \}$. 

If $i\in W_1$  then we choose $j\in W_{2}$  such that $j=\varphi (i)$  and put   $h(x_{2i})=x_{2j},   \   \  $ 
$h(x_{2i+1})=x_{2j+1}. \  $

If $i\in W_{2}$  then we take  $k  = \varphi ^{-1}(i) \ $ and put  $ \ h(x_{2i})= x_{2k}$,   $   \  \ h(x_{2i+1})= x_{2k+1}$.
\vskip 10pt 

By the construction of  $h$, we have $hT_1h=T_2$. Since  $\mathcal{J}$  is translation invariant, we have 
$T_2 \in \mathcal{J}$,  $T \in \mathcal{J}$   so  $\mathcal{J}=\mathcal{I}$ contradicting  $\mathcal{J}\subset \mathcal{I}$. 

To conclude the proof, we put $I =\mathcal{I}^{\wedge}$.  By the construction of  $\mathcal{I}$, $I$ is a closed ideal in  $\beta S_\kappa$ satisfying $(i) $, $(ii) $.   $ \ \  \ \ \Box$
\vskip 15pt 

{\bf Remark 3.2.} 
If $I$  is a subset of  $\beta G$  such that  $\overline{ G ^*G ^*}\subseteq   I $  then $I$  is an ideal in 
$G^*$. It follows that between  $\overline{ G ^*G ^*}$  and  $G^*$  there are no maximal closed ideals in $G^*$. 
\vskip 15pt 

{\bf Lemma 3.1.} {\it
Let  $\{ A_n:  n<\omega \}$  be a family of sparse subsets of a group  $G$,  $A=\cup_{n<\omega} A_n $. Then $A$ is sparse provided that the following two conditions are satisfied :
\vskip 10pt

$(i) $  for every  $F\in \mathcal{F}_G$  there exists  $K\in \mathcal{F}_G$ such that $F(A_i \setminus K)\cap F(A_j \setminus K)= \emptyset $  for all  $i<j<\omega$; 
\vskip 10pt 

$(ii) $ for every  $g\in G \setminus \{e\} $,   there exists $m\in \omega$  such that  $gx\notin A$  for each $x\in \cup _{n>m} A_n$.
\vskip 10pt 

Proof.} We take an arbitrary ultrafilter  $p\in G^*$ and prove that  $\Delta _p (A)$
 is finite. We split the proof in two cases. 
\vskip 10pt 

Case   $\Delta _p (A_n)\neq\emptyset$   for some $n< \omega$. Since $A_n$  is sparse, we have
 $\Delta _p (A_n)= T_p$
  for some  $F\in \mathcal{F}_G$ . We show that $\Delta _p (A)=\Delta _p (A_n)$. Clearly, 
$\Delta _p (A_n)\subseteq  \Delta _p (A)$. 
We take an arbitrary $g\in G\setminus T$,
 put  $F=T\cup \{g\}$
and choose $K$ satisfying $(i) $.  
Then  $A_n\notin gp$ and 
 $\cup \{A_i : i <\omega , \   i\neq n  \}\in gp$ 
so $gp \notin \Delta _p (A)$
and $ \Delta _p (A) \subseteq \Delta _p (A_n)$.
\vskip 10pt 

Case $\Delta _p (A_n)=\emptyset$    for each  $n< \omega$. We show that 
$| \Delta _p (A) | \leq 1$. Assume the contrary : $A\in g_1 p$,   $A\in g_2 p$
for distinct  $g_1 g_2\in G$. We denote  $g=g_1 g_2 ^{-1}$,  $q=g_1 p$. 
Then   $\Delta _p (A)= \Delta _q (A)$, $ \ A\in q$
and  $ \  g^{-1}A\in q$. We choose $m$ satisfying $(ii) $. Since  $\Delta _q (A_n)=\emptyset$  
 for each  $n<\omega$, we have  $\cup_{n>m} (A_n)\in q$
 but $A\notin gq$
 and we get a contradiction with  $ \  g^{-1}A\in q$ .  $ \ \  \ \ \Box$
\vskip 15pt 

{\bf Theorem 3.4.}
{\it  For every group $G$, we have $ {\it FT_G \subset Sp _G}  $ so $\overline{ G ^*G ^*}\subset  { \it FT_G ^{\wedge}  } $.
\vskip 10pt 

Proof. }
Since $ {\it FT_G \subseteq Sp _G} $, 
 we should find a sparse subset $A$ of $G$ which is not $n$-thin for each $n\in  \mathbb N$.  
Passing to a countable subgroup of $G$, we suppose that $G$ itself is countable. 

We construct $A$ in the form  $A=\cup_{n<\omega} A_n$ to satisfy the conditions $(i) $,  $(ii)$  of Lemma 3.1  and such that $A_n$ is not  $n$-thin for each $n>0$.  
For each $n<\omega$, we construct $An$  in the form $A_n = \cup_{i  <\omega}  K_n x_{n_i}$
 for some finite  $K_n$, $ |K_n|= n+1$. $e \in K_n$  and some sequence $(x_{n_i})_{i<\omega}$   in $G$.

We enumerate   $G=\{ g_n : n<\omega \}$,  $g_0=e$   and denote  $F_n = \{ g_n, $ $ \dots ,  g_n \} $.  We put  $K_0=\{ e \}$,     $g_{00}=e$. Assume that we have chosen $K_0, \dots , K_n$
and  $ \{  x_{00}, x_{01},\dots   x_{0n},$   $ \ \dots \ $,     $x_{n0}, x_{n1},\dots   x_{nn} \},$    
so that following  conditions are satisfied: 
\vskip 10pt 

(1)   $ \  \  \{ g_n,   \   \dots ,  \  g_n \} \cap   K_n   K_n^{-1}=\emptyset$;\vskip 7pt 

(2)    $ \  \  F_n  K_m  x_{mn}\cap F_n  K_m $  $\{  x_{m0}, \dots   x_{m \ n-1}\} =\emptyset, \ $   $ \  0\leq m\leq  n$;\vskip 7pt 

(3)   $ \  \  F_n  K_n  \{ x_{n0},  \dots   , x_{n n}\} \cap F_n  K_m $  $\{  x_{m0}, \dots ,  x_{m  n}\} =\emptyset, \ $   $ \  0\leq m<  n$;\vskip 7pt 

(4)   $ \  \  F_n  K_n  x_{ni}\cap F_n  K_n $  $ x_{nj} =\emptyset, \ $   $ \  0\leq i< j\leq  n$.

\vskip 15pt 
Then we choose  $ K_{n+1}$  and $$ \{  x_{0 \ n+1}, x_{ 1 \ n+1},\dots   , x_{n \ n+1} \}, \ \
 \{  x_{ n+1  \ 0}, x_{ n+1  \ 1},\dots   , x_{n+1 \ n+1} \}$$ 
to satisfy (1), (2), (3), (4)  with $n+1$  in place of $n$. After   $\omega$  steps, we  get the family  $\{A
_n : n<\omega \}$.

We put  $K=K_0 \{  x_{00}, x_{01},\dots   x_{0n} \} \bigcup \dots \bigcup K_n
 \{ x_{n0}, x_{n1},\dots  , x_{nn} \}$ . 

By (2), (3), (4), $F_n (A_i \setminus K)\cap F_n (A_j \setminus K)=\emptyset$ for all $i<j<\omega .$  Hence, the condition $(i)$ of Lemma 3.1 is satisfied. 

By (1), (3), (4), $g_ n (\cup_{i>n} A_i) \cap A=\emptyset$, so the condition $(ii)$ is satisfied. By Lemma 3.1, $A$ is sparse.  For every $n<\omega$, the subsets $\{ g \{ x_{n_i} : i< \omega \}: g\in K_n \}$  of $A_n$ are pairwise disjoint. Since $|K_n|=n+1$, $A_n$  is not $n$-thin. $ \ \  \ \ \Box$

For subsets $X, Y$  of a group $G$, we say that the product $XY$ is an $n$-{\it stripe} if $|X|=n$,  
$n\in \mathbb{N}$  and $|X|=\omega$.   It is easy to see that a subset   $A$  of $G$  is $n$-thin  if and only if $A$  has no $(n+1)$-stripes. Thus, $p\in FT_G ^\wedge$ is and only if each member $P\in p$  has an   $n$-stripe for every $n\in \mathbb{N}$. 

We say that $XY$ is an $(n,m)$-{\it rectangle}  if $|X|=n$, $|Y|=m \ $, $ \ n,m\in \mathbb{N}$. We say that a subset $A$ of $G$ {\it has bounded rectangles} if there is $n\in \mathbb{N}$  such that $A$  has no $(n, n)$-rectangles (and so $(n, m)$-rectangles for each $m>n$). 

We denote by $BR_G$ the family of all subsets of $G$ with bounded rectangles. 

\vskip 15pt   

{\bf Theorem 3.5.}
{\it  For a group $G$, the following statements hold:\vskip 10pt 

$(i)$ $BR_G$ is a left translation invariant ideal  in $ \mathcal{P}_G$  ; 

\vskip 10pt 

$(ii)$  $BR_G^\wedge$ is a closed ideal in $\beta G$  and  $p\in BR_G^\wedge$ if and only if each member $P\in p$ has an $(n,n)$-rectangle for  every $n\in   \mathbb{N} $; 

\vskip 10pt 

$(iii)$   $BP_G \subset FT_G$.\vskip 15pt

Proof. } 
$(i) \ $  If $XY$  is an $(n,n)$-rectangle then $(gX)Y$ and $X(Yg)$ are $(n,n)$-rectangles, so the family $BP_G$ is translation invariant.

We take $AB\in BP_G$ and choose $n\in  \mathbb{N}$ such that $A,B$ have no $(n,n)$-rectangles. By the bipartite Ramsey theorem [3, p. 95], there is a natural number $r$ such that, for every 2-coloring of edges of the complete biparte graph $K_{r,r}$, one can find a monochrome copy of  $K_{n,n}$.
 We assume that $A\cup B$  contains an $(r,r)$-rectangle $XY$.  
We define a coloring $\chi : X\times Y \rightarrow \{ 0, 1 \}$  of the Cartesian product $X\times Y$  
by the rule: $\chi ((x, y))=1$ if and only if $xy\in A$. 
By the choice of $r$, there exist $X^\prime\subset X  \  $,  $  \  Y^\prime\subset Y$  such that $| X^\prime |= | Y^\prime |=n \  $ and  $   \  X^\prime \times Y^\prime   \  $ is monochrome. 
Then either $  \  X^\prime Y^\prime\subset A  \  $ or  $  \   X^\prime Y^\prime\subset B  \  $ and we get a contradiction with the choice of $A$ and $B$.  Hence,   $BP_G$  is an ideal in  $\mathcal{P}_G$ .
\vskip 10pt

$(ii \ )$  By $  (i)$,   $ \ BP_G  ^{\wedge}  \  $  is a left ideal and   $ \ (BP_G  ^{\wedge})G \subseteq 
BP_G  ^{\wedge}  $  . 
Since   $ \ BP_G   \subseteq FT_G \subset Sp_G  $  and  $Sp_G^{\wedge}  =  \overline{ G ^*G ^*}  $,
we have  
$ \ (BP_G  ^{\wedge})  G^*   \subseteq BP_G  ^{\wedge}  $  so  $BP_G  ^{\wedge}$ 
is a right ideal. The second statement of $(ii)$ is evident. 
\vskip 10pt

$(iii)$ Passing to subgroups, we suppose that $G$ is countable and construct   $A\in FT_G\setminus  BP_G  $
 in the form   $A=\bigcup _{n<\omega  }   X_N  Y_N  $,   $ \  \  |X_n | = | Y_n| = n+1$. We  enumerate  $G= \{  g_n: n<\omega \}$,    $g_0 =e$
 and  put  $X_0 =  Y_0 = \{0\}$.  Suppose that we have chosen   $X_0   Y_0,  \dots , X_n   Y_n$. We choose 
$X_{n+1}   Y_{n+1}$,   $|X_{n+1} | = | Y_{n+1}| = n+2$
 to satisfy the following conditions for each  $i\in \{1, \dots . n+1 \}$:
$$g_i X_{n+1}  \  Y_{n+1}  \bigcap  X_{n+1}  \  Y_{n+1} = \emptyset ,  \   \    g_i X_{n}  \  Y_{n}  \bigcap  (X_{0}    Y_{0}   \cup  \dots    \cup   X_{n}  \  Y_{n}  )  = \emptyset$$

After $\omega$ steps, we get the desired $A$. Indeed,  $ X_{n}    Y_{n} \subset A$  so  $A\notin BP_G$. By the construction,   $gA\cap  A$  is finite for each  $g\in G\setminus \{e \}$, 
 so $A$  is thin and $A\in   FT_G$.   $ \ \  \ \ \Box$


\section{Between  $\overline{ K(G)}$  and  $\overline{ G ^*G ^*}$ }

Let $(g_n)_{n\in\omega} $ be an injective sequence  in a group $G$. The set 
$$\{   g_{i_1 }  g_{i_2 } \dots g_{i_n } : 0\leq i_1< i_2< \dots < g_{i_n } < \omega \}  $$ is called an {\it $FP$-set }.

Given a sequence $(b_n)_{n\in\omega} $  in $G$, we say that the set 
$$\{   g_{i_1 }  g_{i_2 } \dots g_{i_n }  b_{i_n }  : 0\leq i_1< i_2< \dots < i_n  < \omega \}  $$ is a {\it piecewise shifted  $FP$-set }.

\vskip 10pt

{\bf Theorem 4.1.}
{\it  For a group $G$, the following statements hold:
\vskip 10pt 

$(i)$ $Sc_G^{\wedge} = cl \{\epsilon p:  \epsilon\in G ^*  , \  p\in \beta G,   \   \epsilon\epsilon=\epsilon\}$;

\vskip 10pt 

$(ii)$  $Sc_G^{\wedge}$ is an ideal in $\beta G$ and $p\in Sc_G^{\wedge}$
 if and only if each member of $p$ contains a pierwise shifted  $FP$-set;

\vskip 10pt 

$(iii)$   $Sc_G^{\wedge}$ is  the minimal close ideal in $\beta G$  containing all idempotents of 
$G ^*$.\vskip 15pt

Proof. } 
$(i) \ $
We remind that a subset $A$  of $G$  is scattered if and only if, for each $p\in A^*$, the subset $Gp$ is discrete in $\beta G$.  
Hence, $A$ is not scattered if and only if, there is  $p\in A^*$  such that $Gp$ is not discrete. On the other hand $Gp$ is  not discrete if and only if $p=\epsilon p$ for some idempotent $\epsilon  \in G^*$.\vskip 10pt 

$(ii) \ $    Since $Sc_G$ is a left translation invariant,  $Sc_G ^{\wedge} $ is a left ideal in   $\beta G$. By $(i)  $,   $(Sc_G ^{\wedge})q\subseteq  Sc_G ^{\wedge} $   for each $q\in \beta G$,  so  $Sc_G^{\wedge}$  is a right ideal . \vskip 10pt 

By [1, Theorem 1], a subset $A$ is scattered if and only if $A$ contains no pierwise  shifted $FP$-sets.
 \vskip 10pt 

$(iii) \ $ Let  $\mathcal{M}$  denotes  the minimal closed ideals of   $\beta G$
 containing all idempotents of   $\beta G$. By $(i) $,   $Sc_G^{\wedge} \subseteq   \mathcal{M}$.
Since  $Sc_G^{\wedge} $ is a closed ideal, we have   $\mathcal{M}=Sc_G^{\wedge} $ .  
 $ \   \  \   \    \   \    \Box$
\vskip 15pt

{\bf Remark 4.1.}   If $\mathcal{I}$  is a group ideal in $\mathcal{P}_G$
and  $\mathcal{I}\subseteq   Sp_G  $
 then $\mathcal{I}= \mathcal{F}_G$
(see Remark 3.1).  We can not state  the same if $\mathcal{I}\subseteq   Sc_G $.

Let $G$  be the direct $sum  \  \oplus  _\omega\mathbb{ Z}_{2}$  
of $\omega$  copies of  $\mathbb{ Z}_{2}=\{ 0,1\}$. 
For $g\in G$, we denote by $supt(g)$ the number of non-zero coordinates of $g$. 
We put  $A=\{g\in G: supt(g)=1 \}$
 and consider the minimal group ideal  $\mathcal{I}$   in  $\mathcal{P}_G$ such that $A\in \mathcal{I}$ . If $S\in \mathcal{I}$ then there is $m\in \mathbb{N}$ 
 such that $supt(g) \leq m$ for each  $g\in S$. It follows that $S$ has no piecewise shifted $FP$-sets, so $S$ is scattered and   $\mathcal{I}\subset Sc_G$.
 \vskip 10pt 

The following observation follows directly from the basic properties of multiplication in $\beta G$: each  right shift is continuous and each left shift on element of $g$ is continuous. 
 \vskip 10pt


{\bf Lemma 3.1.} {\it  If $L$ is a left ideal in  $\beta G$ and $R$ is a right ideal in  $\beta G$ then 
$\overline{ LR}$  is an ideal in  $\beta G$.}
\vskip 10pt   

For a group $G$, we put  $I_{G,0} = G^*$ and $I_{G,n+1} = \overline{ G^* I_{G,n }}$.
By  Lemma 4.1, each $I_{G,n}$  is an ideal in   $\beta G$.

Clearly, $I_{G,n+1} \subseteq I_{G,n }$ so $I_{G,n} \subseteq \overline{G^* G^*}$ for $n>0$. 
\vskip 15pt   

{\bf Theorem 4.2.}
{\it  For every group $G$ and $n\in \omega$, we  have \vskip 5pt   

$(i) \ $   $I_{G,n+1} \subset I_{G,n }$ 
\vskip 10pt   

$(ii) \ $   $Sc_{G}^{\wedge} \subset I_{G,n }$. 

\vskip 15pt   

Proof. } 
$(i) \ $  We note that $I_{G,n+1}^{\vee} =\{A\subseteq G:  \Delta _p (A)$ is finite for each $p\in  I_{G,n } \}$ and apply Theorem 4 from [7] stating that  $I_{G,n}^{\vee} \subset I_{G,n+1}^{\vee} $.  
\vskip 10pt   

$(ii) \ $  
For $n=0$, this  is evident.  We take an idempotent $\epsilon \in G^*$,   $p\in  \beta G$ and assume that 
 $\epsilon p\in  I_{G,n-1 }$.
Then   $\epsilon\epsilon p\in  G^* I_{G,n-1 }$,   so  $\epsilon p\in  I_{G,n }$. 
Applying Theorem 4.1, we conclude  that $Se_G ^{\wedge} \subseteq  I_{G,n }$. 
The strict inclusion follows from  $(i)$ .
$ \   \  \   \    \   \    \Box$
\vskip 15pt

For a natural number $n$, we denote by  $(G ^*)^n$
  the product of $n$ copies of $n$. By Lemma 4.1,   $(G ^*)^n$  is an ideal in  $\beta G$.. 
Clearly, $\overline{ (G ^*) ^{n+1}}\subseteq   \overline{ (G ^*) ^{n}}$.
and  $ \overline{ (G ^*) ^{n}}\subseteq I_{G,n }$ .
\vskip 10pt  

By analogy with Theorem 4.2, we can prove \vskip 15pt  

{\bf Theorem 4.3.}
{\it  For every  group $G$  and $n\in\omega$, we have   \vskip 5pt 

$(i) \ $   $\overline{ (G ^*) ^{n+1}} \subset  \overline{ (G ^*) ^{n}}$;
\vskip 5pt   

$(ii) \ $   $Sc_{G}^{\wedge} \subset  \overline{ (G ^*) ^{n}}$.} 

\vspace{3 mm}

\centerline{\bf References }

\vspace{3 mm}

[1]  T. Banakh, I. Protasov, S. Slobodianiuk, {\it Scattered subsets of groups}, Ukr. Math. J. {\bf 67} (2015),  347-356. 

[2]	{M. Filali, Ie. Lutsenko, I. Protasov,} {\em Boolean group ideal and the ideal structure of $\beta G$},  Mat. Stud.
{\bf 31} (2009), 19-28.
  
[3] {R. Graham, B. Rotschild, J. Spencer,}  {\em Ramsey  Theory}, Willey, New York,  
1980.

[4] { N. Hindman, D. Strauss} {\em  Algebra in the Stone-$\check{C}$ech compactification},   de Gruyter, Berlin, New York, 1998.

[5] {Ie. Lutsenko, I. Protasov,} {\em   Sparse, thin and other subsets of groups},  Algebra Computa-tion {\bf 19} (2009), 491-510.

[6]	{ Ie. Lutsenko, I. Protasov,} {\em Thin subsets of balleans},  Appl. Gen. Topology
{\bf 11} (2010), 89-93.

[7]	{ Ie. Lutsenko, I. Protasov,} {\em Relatively thin and sparse subsets of  groups},  Ukr. Math. J. {\bf 63} (2011), 254-264.

[8]	{ O. Petrenko, I. Protasov,} {\em Thin ultrafilters},  Notre Dame J. Formal Logic
{\bf 53} (2012), 79-88.

[9] {I. Protasov,}  {\em Selective survey  on subset combinatories of  groups, } J. Math. Sciences  
{\bf  174} (2011), 486-514.

[10] I. Protasov, T. Banakh, {\it Ball structures and colorings of groups and graphs},  Math. Stud. Monogr. Ser., Vol. 11, VNTL, Lviv, 2003.

[11] {I. Protasov,  K. Protasova, }  {\em Ramsey-product subsets of a group, },  preprint (arxiv: 1703, 03874).

[12] {I. Protasov,  O. Protasova, }  {\em Sketch of group balleans, } Math. Stud. {\bf  22} (2004), 10-20.

[13] {I. Protasov,  O. Protasova, }  {\em On closed ideals in $\beta G$}, Semigroup Forum  {\bf  75} (2007), 237-240.

[14] I. Protasov, S. Slobodianiuk, {\it  Thin subsets of  groups}, Ukr.  Math. J. {\bf 65} (2013), 1384-1393.

[15]  I. Protasov, S. Slobodianiuk, {\it On the subset combinatorics of $G$-spaces},  Algebra Discrete Math.  {\bf 17} (2014), No 1, 98-109.

[16] {I. Protasov,  S. Slobodianiuk, }  {\em Ultracompanions of subsets of a group, } Comment. Math. Univ. Carolin. {\bf  55} (2014), 257-265.

\vspace{5 mm}

CONTACT INFORMATION

I.~Protasov: \\
Faculty of Computer Science and Cybernetics  \\
        Kyiv University  \\
         Academic Glushkov pr. 4d  \\
         03680 Kyiv, Ukraine \\ i.v.protasov@gmail.com

\medskip

K.~Protasova:\\
Faculty of Computer Science and Cybernetics \\
        Kyiv University  \\
         Academic Glushkov pr. 4d  \\
         03680 Kyiv, Ukraine \\ ksuha@freenet.com.ua

\end{document}